\documentclass[a4paper,12pt]{amsart}
\usepackage{amssymb}
\usepackage{ifthen}
 \usepackage[dvips]{graphicx}
 \usepackage{blindtext}
 \usepackage{mathtools}
\nonstopmode \numberwithin{equation}{section}
\usepackage{amssymb}
\usepackage{ifthen}
\usepackage{graphicx}
\usepackage{amsmath}
\usepackage[T1]{fontenc} %skandit
\usepackage[utf8]{inputenc}
\usepackage[usenames,dvipsnames]{color}
\usepackage{color}
\usepackage[english]{babel}
\usepackage{fancyhdr}
\usepackage{fancybox}
\usepackage{tikz}

\nonstopmode \numberwithin{equation}{section}
\theoremstyle{plain}
\newtheorem{thm}[equation]{Theorem}
\newtheorem{cor}[equation]{Corollary}
\newtheorem{lem}[equation]{Lemma}
\newtheorem{prop}{Proposition}

\newtheorem{conj}{Conjecture}

\theoremstyle{definition}
\newtheorem{defn}{Definition}[section]

\newtheorem{prob}{Problem}
\newtheorem{rem}{Remark}[section]

\newcounter{minutes}
\divide\time by 60
\newcounter{hours}
\multiply\time by 60
\addtocounter{minutes}{-\time}
\newcounter {own}
\def\theown {\thesection       .\arabic{own}}

\newenvironment{pf}[1][]{%
 \vskip 3mm
 \noindent
 \ifthenelse{\equal{#1}{}}%
  {{\slshape Proof. }}%
  {{\slshape #1.} }%
 }%
{\qed\bigskip}

\newcounter{alphabet}

\newcommand{\real}{{\operatorname{Re}\,}}

\def\be{\begin{equation}}
\def\ee{\end{equation}}

\newcommand{\bee}{\begin{enumerate}}
\newcommand{\eee}{\end{enumerate}}

\newcommand{\blem}{\begin{lem}}
\newcommand{\elem}{\end{lem}}
\newcommand{\bthm}{\begin{thm}}
\newcommand{\ethm}{\end{thm}}
\newcommand{\bcor}{\begin{cor}}
\newcommand{\ecor}{\end{cor}}
\newcommand{\beg}{\begin{examp}}
\newcommand{\eeg}{\end{examp}}
\newcommand{\begs}{\begin{examples}}
\newcommand{\eegs}{\end{examples}}

\newcommand{\bdefn}{\begin{defn}}
\newcommand{\edefn}{\end{defn}}

\newcommand{\bprob}{\begin{prob}}
\newcommand{\eprob}{\end{prob}}
\newcommand{\bei}{\begin{itemize}}
\newcommand{\eei}{\end{itemize}}

\newcommand{\bcon}{\begin{conj}}
\newcommand{\econ}{\end{conj}}
\newcommand{\bcons}{\begin{conjs}}
\newcommand{\econs}{\end{conjs}}
\newcommand{\bprop}{\begin{prop}}
\newcommand{\eprop}{\end{prop}}
\newcommand{\br}{\begin{rem}}
\newcommand{\er}{\end{rem}}
\newcommand{\brs}{\begin{rems}}
\newcommand{\ers}{\end{rems}}
\newcommand{\bo}{\begin{obser}}
\newcommand{\eo}{\end{obser}}
\newcommand{\bos}{\begin{obsers}}
\newcommand{\eos}{\end{obsers}}
\newcommand{\bpf}{\begin{pf}}
\newcommand{\epf}{\end{pf}}
\newcommand{\ba}{\begin{array}}
\newcommand{\ea}{\end{array}}
\newcommand{\beq}{\begin{eqnarray}}
\newcommand{\beqq}{\begin{eqnarray*}}
\newcommand{\eeq}{\end{eqnarray}}
\newcommand{\eeqq}{\end{eqnarray*}}

\begin{document}

\title{Bohr radius for certain close-to-convex harmonic mappings}

\author{Molla Basir Ahamed}
\address{Molla Basir Ahamed,
	School of Basic Science,
	Indian Institute of Technology Bhubaneswar,
	Bhubaneswar-752050, Odisha, India.}
\email{mba15@iitbbs.ac.in}

\author{Vasudevarao Allu}
\address{Vasudevarao Allu,
School of Basic Science,
Indian Institute of Technology Bhubaneswar,
Bhubaneswar-752050, Odisha, India.}
\email{avrao@iitbbs.ac.in}

\author{Himadri Halder}
\address{Himadri Halder,
School of Basic Science,
Indian Institute of Technology Bhubaneswar,
Bhubaneswar-752050, Odisha, India.}
\email{hh11@iitbbs.ac.in}

\subjclass[{AMS} Subject Classification:]{Primary 30C45, 30C50, 30C80}

\def\thefootnote{}
\footnotetext{ {\tiny File:~\jobname.tex,
printed: \number\year-\number\month-\number\day,
          \thehours.\ifnum\theminutes<10{0}\fi\theminutes }
} \makeatletter\def\thefootnote{\@arabic\c@footnote}\makeatother

\begin{abstract}
 Let $ \mathcal{H} $ be the class of harmonic functions $ f=h+\bar{g} $ in the unit disk $\mathbb{D}:=\{z\in\mathbb{C} : |z|<1\}$, where $ h $ and $ g $ are analytic in $ \mathbb{D} $. Let 
 $$\mathcal{P}_{\mathcal{H}}^{0}(\alpha)=\{f=h+\overline{g} \in \mathcal{H} : \real (h^{\prime}(z)-\alpha)>|g^{\prime}(z)|\; \mbox{with}\; 0\leq\alpha<1,\; g^{\prime}(0)=0,\;  z \in \mathbb{D}\}
 $$ be the class of close-to-convex mappings defined by Li and Ponnusamy \cite{Injectivity section}. In this paper, we obtain the sharp Bohr-Rogosinski radius, improved Bohr radius and refined Bohr radius for the class $ \mathcal{P}_{\mathcal{H}}^{0}(\alpha) $.
\end{abstract}

\maketitle
\pagestyle{myheadings}
\markboth{Molla Basir Ahamed, Vasudevarao Allu and  Himadri Halder}{Bohr radius for certain close-to-convex harmonic mappings}

\section{Introduction}
The classical inequality of Bohr says that if $f$ is an analytic function in the unit disk $\mathbb{D}:=\{z\in \mathbb{C}: |z|<1\}$ with the following Taylor series expansion
\begin{equation} \label{him-p4-e-1.1}
	f(z)=\sum_{n=0}^{\infty} a_{n}z^{n}
\end{equation}
such that $|f(z)|<1$ in $\mathbb{D}$, then the majorant series $M_{f}(r)$ associated with $f$ satisfies the following inequality
\begin{equation} \label{him-p4-e-1.2}
	M_{f}(r):= \sum_{n=0}^{\infty} |a_{n}|r^{n} \leq 1 \quad \mbox{for} \quad |z|=r \leq 1/3,
\end{equation}
and the constant $1/3$, known as Bohr radius, cannot be improved. In 1914, H. Bohr \cite{Bohr-1914} obtained the inequality \eqref{him-p4-e-1.2} for $r \leq 1/6$ and subsequently improved by $ 1/3 $ later, Weiner, Riesz and Schur independently obtained the constant $1/3$. An observation shows that the quantity $1-|a_{0}|$ is equal to $d(f(0),\partial f(\mathbb{D}))$. Therefore, the inequality \eqref{him-p4-e-1.2} is called Bohr inequality, can be written in the following form 
\begin{equation} \label{him-p4-e-1.3}
	\sum_{n=1}^{\infty} |a_{n}z^{n}|\leq 1-|a_{0}|=d(f(0),\partial f(\mathbb{D}))
\end{equation}
for $|z|=r \leq 1/3$, where $d$ is the Euclidean distance. It is important to note that the constant $1/3$ is independent of the coefficients of the Taylor series \eqref{him-p4-e-1.1}. This fact can be elucidated by saying that Bohr inequality occurs in the class $\mathcal{B}$ of analytic self maps of the unit disk $\mathbb{D}$. Analytic functions $f \in \mathcal{B}$ of the form \eqref{him-p4-e-1.1} satisfying the inequality \eqref{him-p4-e-1.2} for $|z|=r\leq1/3$, are sometimes  said to satisfy the classical Bohr phenomenon. The notion of the Bohr phenomenon can be generalized to the class $\mathcal{F}$ consisting of analytic functions $f$ from $\mathbb{D}$ to a given domain $\varOmega \subseteq \mathbb{C}$ such that $f(\mathbb{D}) \subseteq \varOmega$ and the class $\mathcal{F}$ is said to satisfy the Bohr phenomenon if there exists largest radius $r_{\varOmega} \in (0,1)$ such that the inequality \eqref{him-p4-e-1.3} holds for $|z|=r\leq r_{\varOmega}$ and for all functions $f\in \mathcal{F}$. We say the largest radius $r_{\varOmega}$ is the Bohr radius for the class $\mathcal{F}$. The Bohr radius has been obtained for the class $\mathcal{F}$ when $\varOmega$ is convex domain \cite{aizn-2007}, simply connected domain \cite{Abu-2010}, the exterior of the closed unit disk, the punctured unit disk, and concave wedge domain (see \cite{Ali-2017}). In $ 1997 $, Boas and Khavinson \cite{Boas-Khavinson-PAMS-1997} generalized the Bohr inequality in several complex variables by finding multidimensional Bohr radius. In $ 2020 $, Liu and Ponnusamy \cite{Liu-Ponnusamy-PAMS-2020} obtained multidimensional analogues of refined Bohr inequality.

\vspace{4mm}
There are many improved versions of Bohr's inequality \eqref{him-p4-e-1.2} in various forms obtained by several authors. In $ 2020 $, Kayumov and Ponnusamy \cite{kayumov-2018-c} obtained several interesting improved versions of Bohr inequality. For more results on this, we refer the reader to glance through the articles  (see \cite{Huang-Ponnusamy-2020,Ismagilov-2020,kayumov-2017,kayumov-2019,kayumov-2018-c,Liu-2020,Ponnusamy-Vijaya-ResultsMath-2020,Ponnusamy-CMFT-2020}). In 2017, Kayumov and Ponnusamy \cite{kayumov-2017} introduced Bohr-Rogosinski radius motivated by Rogosinski radius for bounded analytic functions in $\mathbb{D}$. Rogosinski radius is defined as follows: Let $f(z)=\sum_{n=0}^{\infty} a_{n}z^{n}$ be analytic in $\mathbb{D}$ and its corresponding partial sum of $f$ is defined by $S_{N}(z):=\sum_{n=0}^{N-1} a_{n}z^{n}$. Then, for every $N \geq 1$, we have $|\sum_{n=0}^{N-1} a_{n}z^{n}|<1$ in the disk $|z|<1/2$ and the radius $1/2$ is sharp. Motivated by Rogosinski radius, Kayumov and Ponnusamy have considered the Bohr-Rogosinski sum $R_{N}^{f}(z)$ which is defined by 
\begin{equation}
	R_{N}^{f}(z):=|f(z)|+ \sum_{n=N}^{\infty} |a_{n}||z|^{n}.
\end{equation}
It is worth to point out that $|S_{N}(z)|=\big|f(z)-\sum_{n=N}^{\infty} a_{n}z^{n}\big| \leq |R_{N}^{f}(z)|$. Therefore, it is easy to see that the validity of Bohr-type radius for $R_{N}^{f}(z)$, which is related to the classical Bohr sum (Majorant series) in which $f(0)$ is replaced by $f(z)$, gives Rogosinski radius in the case of bounded analytic functions in $\mathbb{D}$. We have the following interesting results by Kayumov and Ponnusamy \cite{kayumov-2017}.

\begin{thm} \cite{kayumov-2017} \label{him-thm-1.5}
	Let $f(z)=\sum_{n=0}^{\infty} a_{n}z^{n}$ be analytic in $\mathbb{D}$ and $|f(z)|\leq 1$. Then
	\begin{equation} \label{him-e-1.6}
		|f(z)|+\sum_{n=N}^{\infty}|a_n||z|^n\leq1
	\end{equation}
	for $|z|=r \leq R_{N}$, where $R_{N}$ is the positive root of the equation $\psi _{N}(r)=0$, $\psi _{N}(r)=2 (1+r)r^{N}-(1-r)^{2}$. The radius $R_{N}$ is the best possible. Moreover, 
	\begin{equation} \label{him-e-1.7}
		|f(z)|^{2}+\sum_{n=N}^{\infty}(|a_n|+|b_n|)|z|^n\leq1
	\end{equation}
	for $R'_{N}$, where $R'_{N}$ is the positive root of the equation $ (1+r)r^{N}-(1-r)^{2}$. The radius $R'_{N}$ is the best possible.
\end{thm}
Recently, Kayumov and Ponnusamy \cite{kayumov-2017} have proved the following improved version of Bohr's inequality.
\begin{thm}\cite{kayumov-2017}
	Let $f(z)=\sum_{n=0}^{\infty} a_{n}z^{n}$ be analytic in $\mathbb{D}$, $|f(z)|\leq 1$ and $S_{r}$ denote the image of the subdisk $|z|<r$ under mapping $f$. Then 
	\begin{equation}
		B_{1}(r):=\sum_{n=0}^{\infty}|a_n|r^n+ \frac{16}{9} \left(\frac{S_{r}}{\pi}\right) \leq 1 \quad \mbox{for} \quad r \leq \frac{1}{3}
	\end{equation}
	and the numbers $1/3$, $16/9$ cannot be improved. Moreover, 
	\begin{equation}
		B_{2}(r):=|a_{0}|^{2}+\sum_{n=1}^{\infty}|a_n|r^n+ \frac{9}{8} \left(\frac{S_{r}}{\pi}\right) \leq 1 \quad \mbox{for} \quad r \leq \frac{1}{2}
	\end{equation}
	and the numbers $1/2$, $9/8$ cannot be improved.
\end{thm}

\vspace{4mm}
Bohr's phenomenon for the complex-valued harmonic mappings have been studied extensivel by many authors (see \cite{Abu-2010,abu-2014,Himadri-Vasu-P1,Himadri-Vasu-P4}). Improved Bohr inequality for locally univalent harmonic mappings have been discussed by Evdoridis \textit{et al.} \cite{evdoridis-2018}. 
\vspace{2mm}
\par A complex-valued function $f=u+iv$ is harmonic if $u$ and $v$ are real-harmonic in $\mathbb{D}$. Every harmonic function $f$ has the canonical representation $f=h+ \overline{g}$, where $h$ and $g$ are analytic in $\mathbb{D}$ known respectively as the analytic and co-analytic parts of $ f $. A locally univalent harmonic function $f$ is said to be sense-preserving if the Jacobian of $f$, defined by $J_{f}(z):=|h'(z)|^{2}-|g'(z)|^{2}$, is positive in $\mathbb{D}$ and sense-reversing if $J_{f}(z)$ is negative in $\mathbb{D}$. Let $\mathcal{H}$ be the class of all complex-valued harmonic functions $f=h+\overline{g}$ defined in $\mathbb{D}$, where $h$ and $g$ are analytic 
in $\mathbb{D}$ such that $h(0)=h'(0)-1=0$ and $g(0)=0$. A function $f \in \mathcal{H}$ is said to be in $\mathcal{H}_{0}$ if $g'(0)=0$.
Thus, every $f=h+\overline{g}\in \mathcal{H}_{0}$ has the following form 
\begin{equation}\label{e-1.11a}
	f(z)=h(z)+\overline{g(z)}=z+\sum_{n=2}^{\infty}a_nz^n+\overline{\sum_{n=2}^{\infty}b_nz^n}.
\end{equation}
In 2013, Ponnusamy {\it et al.} \cite{S.Ponnusamy variability regions-2013}  considered 
the following classes   
$$\mathcal{P}_{\mathcal{H}}^{0}=\{f=h+\bar{g} \in \mathcal{H} : \real h'(z)>|g'(z)| \quad \mbox{with} \quad g'(0)=0 \quad  \mbox{for} \quad z \in \mathbb{D}\}.
$$
Motivated by the class $\mathcal{P}_{\mathcal{H}}^{0}$, Li and Ponnusamy \cite{Injectivity section} have studied the following class $\mathcal{P}_{\mathcal{H}}^{0}(\alpha)$ defined by
$$\mathcal{P}_{\mathcal{H}}^{0}(\alpha)=\{f=h+\overline{g} \in \mathcal{H} : \real (h^{\prime}(z)-\alpha)>|g^{\prime}(z)|\; \mbox{with}\; 0\leq\alpha<1,\; g^{\prime}(0)=0\; \mbox{for}\;  z \in \mathbb{D}\}.
$$
We have the following coefficient bounds and growth estimates for the class $\mathcal{P}_{H}^{0}(\alpha)$.
\begin{lem} \label{lem-1.8}  \cite{Injectivity section}
	Let $f \in \mathcal{P}^{0}_{\mathcal{H}}(\alpha) $ and be given by \eqref{e-1.11a}. Then for any $n \geq 2$, 
	\begin{enumerate}
		\item[(i)] $\displaystyle |a_n| + |b_n|\leq \frac {2(1-\alpha)}{n}; $\\[2mm]
		
		\item[(ii)] $\displaystyle ||a_n| - |b_n||\leq \frac {2(1-\alpha)}{n};$\\[2mm]
		
		\item[(iii)] $\displaystyle |a_n|\leq \frac {2(1-\alpha)}{n}.$
	\end{enumerate}
	All the inequalities  are sharp, with extremal function $f(z)=(1-\alpha)(-z-2\,\;\log(1-z))+\alpha z$.
\end{lem}
\begin{lem} \cite{Himadri-Vasu-P1}\label{lem-1.9}
	Let $f=h+\overline{g} \in \mathcal{P}^{0}_{\mathcal{H}}(\alpha)$ with $0\leq \alpha <1$. Then 
	\begin{equation} \label{e-1.10}
		|z|+ \sum\limits_{n=2}^{\infty}  \dfrac{2(1-\alpha)(-1)^{n-1}}{n} |z|^{n} \leq |f(z)| \leq |z|+ \sum\limits_{n=2}^{\infty}  \dfrac{2(1-\alpha)}{n} |z|^{n}.
	\end{equation} 
	Both  inequalities are sharp.
\end{lem}
The organization of this paper is follows: In section 2 we obtain sharp Bohr-Rogosinski radius for the class $ \mathcal{P}^{0}_{\mathcal{H}}(\alpha) $ of close-to-convex harmonic mappings. In section 3, we establish interesting sharp improved-Bohr radius $ \mathcal{P}^{0}_{\mathcal{H}}(\alpha) $. In section 4, we prove sharp refined-Bohr radius as well as Bohr-type inequality for the class $ \mathcal{P}^{0}_{\mathcal{H}}(\alpha) $. Section 6 is devoted for the proofs of main the results.
\section{Bohr-Rogosinski Radius for the class $ \mathcal{P}^{0}_{\mathcal{H}}(\alpha) $ }
We first prove the following Bohr-Rogosinski radius for the class $ \mathcal{P}^{0}_{\mathcal{H}}(\alpha) $.
\begin{thm}\label{th-2.1}
	Let $ f\in \mathcal{P}^{0}_{\mathcal{H}}(\alpha) $ be given by \eqref{e-1.11a}. Then, for $ N\geq 2 $,
	\begin{equation}\label{e-2.2}
		|f(z)|+\sum_{n=N}^{\infty}(|a_n|+|b_n|)|z|^n\leq d\left(f(0),\partial f(\mathbb{D})\right)
	\end{equation} 
	for $ |z|=r\leq r_{_N}(\alpha) $, where $ r_{_N}(\alpha) $ is the smallest root of the equation
	\begin{equation}\label{e-2.3}
		r-1-2(1-\alpha)\left(r-1+\ln\,(2(1-r)^2)+\sum_{n=1}^{N-1}\frac{r^n}{n}\right)=0\;\; \mbox{in}\;\; (0,1).
	\end{equation}
The constant $ r_{_N}(\alpha) $ is the best possible.
\end{thm}
\vspace{0.5in}
\begin{figure}[!htb]
	\begin{center}
		\includegraphics[width=0.40\linewidth]{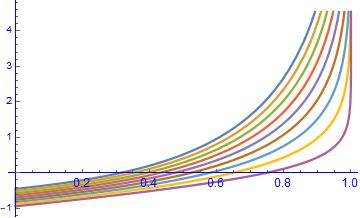}\;\;\;\;\;\;\;\;\;\includegraphics[width=0.40\linewidth]{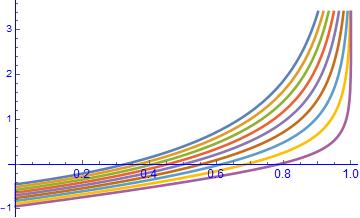}
	\end{center}
	\caption{The graph  of $ r_3(\alpha) $ and $ r_{10}(\alpha) $ in $ (0,1) $.}
\end{figure}

\begin{table}[ht]
	\centering
	\begin{tabular}{|l|l|l|l|l|l|l|l|l|l|}
		\hline
		$\;\;\alpha$& $\;\;0.1$&$\;\;0.2$& $\;\;0.3$& $\;\; 0.4 $& $\;\;0.5$&$\;\;0.6$& $\;\;0.7$& $ \;\;0.8 $& $ \;\;0.9 $ \\
		\hline
		$r_2(\alpha)$& $0.2771 $&$0.3115 $& $0.3477$& $0.3866$& $0.4296 $&$0.4785 $& $0.5367$& $0.6109$ &$ 0.7187 $\\
		\hline
		$r_3(\alpha)$& $0.3121 $&$0.3493 $& $0.3877$& $0.4281$& $0.4717 $&$0.5201 $& $0.5764$& $0.6463$ &$ 0.7453 $\\
		\hline
		$r_{6}(\alpha)$& $0.3248 $&$0.3653 $& $0.4070$& $0.4508$& $0.4978 $&$0.5493 $& $0.6080$& $0.6786$ &$ 0.7736 $\\
		\hline
		$r_{10}(\alpha)$& $0.3251 $&$0.3657 $& $0.4978$& $0.4522$& $0.4999 $&$0.5527 $& $0.6130$& $0.6859$ &$ 0.7832 $\\
		\hline
	\end{tabular}
	\vspace{1mm}
	\caption{This table shows the value of the roots $ r_{_N}(\alpha) $ for different values of $ \alpha $ in $ [0,1) $ and $ N=2, 3, 6, 10 $.}
	\label{tabel-1}
\end{table}
\begin{thm}\label{th-2.4}
	Let $ f\in \mathcal{P}^{0}_{\mathcal{H}}(\alpha) $ be given by \eqref{e-1.11a}. Then, $ N\geq 2 $,
	\begin{equation}\label{e-2.5}
		|f(z)|^2+\sum_{n=N}^{\infty}|a_n||z|^n\leq d\left(f(0),\partial f(\mathbb{D})\right)
	\end{equation} 
	for $ |z|=r\leq r_{_N}(\alpha) $, where $ r_{_N}(\alpha)\in (0,1) $ is the smallest root of the equation
	\begin{equation}\label{e-2.6}
		\bigg(r-2(1-\alpha)(r+\ln\;(1-r))\bigg)^2-2(1-\alpha)\left(\ln\,(2-2r)-1+\sum_{n=2}^{N-1}\frac{r^n}{n}\right)-1=0.
	\end{equation}
The constant $ r_{_N}(\alpha) $ is the best possible.
\end{thm}

\begin{figure}[!htb]
	\begin{center}
		\includegraphics[width=0.40\linewidth]{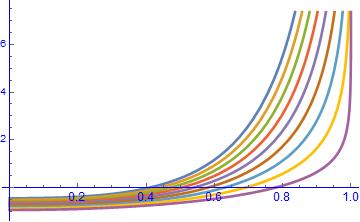}\;\;\;\;\;\;\;\;\;\;\;\;\;\includegraphics[width=0.40\linewidth]{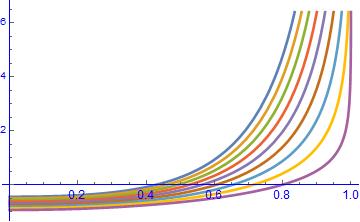}
	\end{center}
	\caption{The graph of $ r_3(\alpha) $ and $ r_8(\alpha) $ of the equation \eqref{e-2.6}.}
\end{figure}
\begin{table}[ht]
	\centering
	\begin{tabular}{|l|l|l|l|l|l|l|l|l|l|}
		\hline
		$\;\;\alpha$& $\;\;0.1$&$\;\;0.2$& $\;\;0.3$& $\;\; 0.4 $& $\;\;0.5$&$\;\;0.6$& $\;\;0.7$& $ \;\;0.8 $& $ \;\;0.9 $ \\
		\hline
		$r_3(\alpha)$& $0.4102 $&$0.4399 $& $0.4708$& $0.5038$& $0.5399 $&$0.5807 $& $0.6291$& $0.6903$ &$ 0.7783 $\\
		\hline
		$r_8(\alpha)$& $0.4304 $&$0.4613 $& $0.4933$& $0.5273$& $0.5644 $&$0.6060 $& $0.6547$& $0.7152$ &$ 0.7994 $\\
		\hline
	\end{tabular}
\vspace{0.08in}
	\caption{Values of $ r_N(\alpha) $ for $ N=3 $ and $ 8 $ when $ \alpha\in[0,1) $.}
	\label{tabel-2}
\end{table}
\begin{thm}\label{th-2.7}
	Let $ f\in \mathcal{P}^{0}_{\mathcal{H}}(\alpha) $ be given by \eqref{e-1.11a}. Then for a positive integer $ N\geq 2 $,
	\begin{equation}\label{e-2.8}
		|f(z^m)|+\sum_{n=N}^{\infty}|a_n||z|^n\leq d\left(f(0),\partial f(\mathbb{D})\right)
	\end{equation} 
	for $ |z|=r\leq r_{m,N}(\alpha) $, where $ r_{m,N}(\alpha)\in (0,1) $ is the smallest root of the equation
	\begin{equation}\label{e-2.9}
		r^m-1-2(1-\alpha)\left(r^m-1+\ln\,((1-r^m)(2-2r))+\sum_{n=1}^{N-1}\frac{r^n}{n}\right)=0.
	\end{equation}
The constant $ r_{m,N}(\alpha) $ is the best possible.
\end{thm}

\begin{figure}[!htb]
	\begin{center}
		\includegraphics[width=0.30\linewidth]{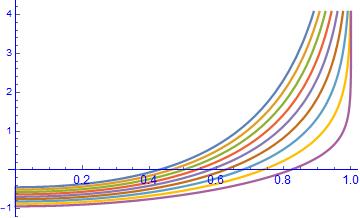}\;\;\;\;\;\includegraphics[width=0.30\linewidth]{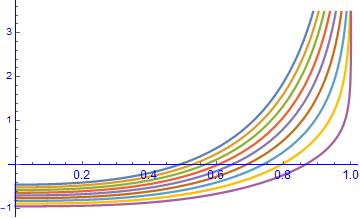}\;\;\;\;\;\includegraphics[width=0.30\linewidth]{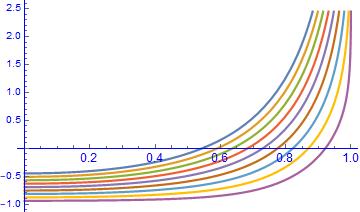}
	\end{center}
\caption{The graphs $ r_{2,2}(\alpha) $, $ r_{3,2}(\alpha) $ and  $ r_{7,2}(\alpha) $ of the equation \eqref{e-2.9}}
	\begin{center}
		\includegraphics[width=0.30\linewidth]{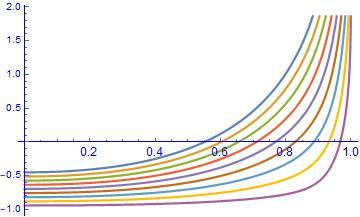}\;\;\;\;\;\includegraphics[width=0.30\linewidth]{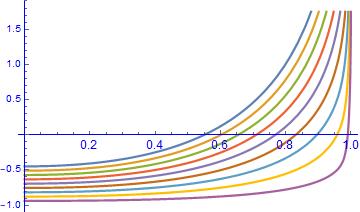}\;\;\;\;\;\includegraphics[width=0.30\linewidth]{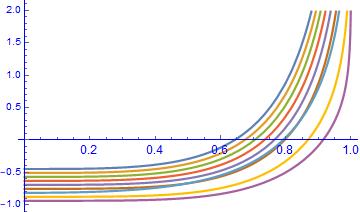}
	\end{center}
\caption{The graphs $ r_{25,2}(\alpha) $, $ r_{150,2}(\alpha) $ and $ r_{5,3}(\alpha) $ of the equation \eqref{e-2.9}.}
%\vspace{0.1in}
	\begin{center}
		\includegraphics[width=0.30\linewidth]{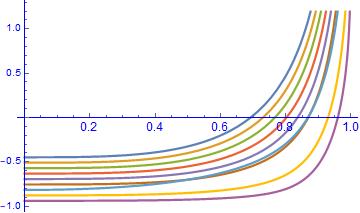}\;\;\;\;\;\includegraphics[width=0.30\linewidth]{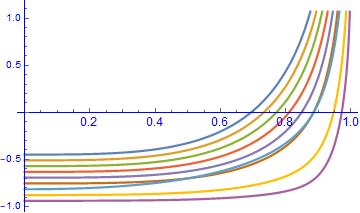}\;\;\;\;\;\includegraphics[width=0.30\linewidth]{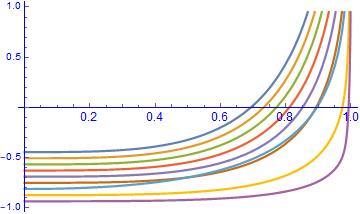}
	\end{center}
\caption{The graphs $ r_{15,3}(\alpha) $, $ r_{25,3}(\alpha) $ and $ r_{180,3}(\alpha) $ of the equation \eqref{e-2.9}.}
%\vspace{0.1in}
	\begin{center}
		\includegraphics[width=0.30\linewidth]{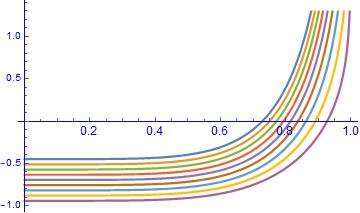}\;\;\;\;\;\includegraphics[width=0.30\linewidth]{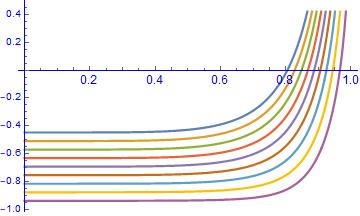}\;\;\;\;\;\includegraphics[width=0.30\linewidth]{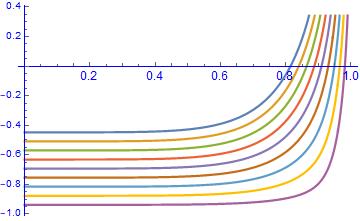}
	\end{center}
\caption{The graphs $ r_{5,5}(\alpha) $, $ r_{15,5}(\alpha) $ and $ r_{35,5}(\alpha) $ of the equation \eqref{e-2.9}.}

\end{figure}

\begin{table}[ht]
	\centering
	\begin{tabular}{|l|l|l|l|l|l|l|l|}
		\hline
		$\alpha$& $r_{2,2}(\alpha)$&$r_{3,2}(\alpha)$& $r_{7,2}(\alpha)$& $ r_{25,2}(\alpha) $ & $ r_{50,2}(\alpha) $& $ r_{90,2}(\alpha) $& $ r_{150,2}(\alpha) $ \\
		\hline
		$0.1$& $0.2016 $&$0.2157 $& $0.2201$& $0.2201$& $0.2201$& $0.2202$ & $0.2202$\\
		\hline
		$0.2$& $0.2436 $&$0.2639$& $0.2724$& $0.2725$& $0.2725$& $0.2725$& $0.2725$\\
		\hline
		$0.3$& $0.2905 $&$0.3187$& $0.3344$& $0.3346$& $0.3346$& $0.3346$& $0.3346$\\
		\hline
		$0.4$& $0.3433 $&$0.3805$& $0.4083$& $0.4093$& $0.4093$& $0.4093$& $0.4093$\\
		\hline
		$0.5$& $0.4030 $&$0.4499$& $0.4962$& $0.5000$& $0.5000$& $0.5000$& $0.5000$\\
		\hline
		$0.6$& $0.4710 $&$0.5268$& $0.5970$& $0.6105$& $0.6106$& $0.6106$& $0.6106$\\
		\hline
		$0.7$& $0.5498 $&$0.6119$& $0.7031$& $0.7430$& $0.7433$& $0.7433$& $0.7433$\\
		\hline
		$0.8$& $0.6443 $&$0.7069$& $0.8033$& $0.8772$& $0.8877$& $0.8884$& $0.8884$\\
		\hline
		$0.9$& $0.7667 $&$0.8184$& $0.8944$& $0.9542 $& $0.9707$& $0.9796$& $0.9848$\\
		\hline
	\end{tabular}
	\vspace{1.5mm}
	\caption{The roots of $ r_{m,N}(\alpha) $ for different values of $ m=2, 3, 7, 25, 50, 90, 150 $ with $ N=2 $ and $ \alpha\in [0,1) $.}
	\label{tabel-4}
\end{table}
\begin{table}[ht]
	\centering
	\begin{tabular}{|l|l|l|l|l|l|l|l|l|}
		\hline
		$\alpha$& $r_{5,3}(\alpha)$&$r_{15,3}(\alpha)$& $r_{35,3}(\alpha)$& $ r_{85,3}(\alpha) $ & $ r_{180,3}(\alpha) $& $ r_{5,5}(\alpha) $& $ r_{15,5}(\alpha) $& $ r_{35,5}(\alpha) $ \\
		\hline
		$0.1$& $0.6435 $&$0.6922 $& $0.6936$& $0.6936$& $0.6936$& $0.7283$ & $0.8048$& $0.8145$\\
		\hline
		$0.2$& $0.6744 $&$0.7298$& $0.7326$& $0.7326$& $0.7326$& $0.7503$& $0.8371$& $0.8397$\\
		\hline
		$0.3$& $0.7045 $&$0.7664$& $0.7716$& $0.7716$& $0.7716$& $0.7717$& $0.8479$& $0.8639$\\
		\hline
		$0.4$& $0.7344 $&$0.8019$& $0.8109$& $0.8111$& $0.8111$& $0.7930$& $0.8678$& $0.8872$\\
		\hline
		$0.5$& $0.7647 $&$0.8360$& $0.8508$& $0.8515$& $0.8515$& $0.8147$& $0.8870$& $0.9094$ \\
		\hline
		$0.6$& $0.7902 $&$0.8687$& $0.8901$& $0.8930$& $0.8930$& $0.8374$& $0.9057$& $0.9300$\\
		\hline
		$0.7$& $0.8299 $&$0.8996$& $0.9257$& $0.9345$& $0.9349$& $0.8621$& $0.9345$& $0.9486$\\
		\hline
		$0.8$& $0.8677 $&$0.9293$& $0.9545$& $0.9686$& $0.9732$& $0.8905$& $0.9440$& $0.9653$ \\
		\hline
		$0.9$& $0.9141 $&$0.9593$& $0.9771$& $0.9875 $& $0.9924$& $0.9267$& $0.9657$& $0.9807$\\
		\hline
	\end{tabular}
	\vspace{0.5mm}
	\caption{This table shows the roots $ r_{m,N}(\alpha) $ for different values of $ m=5, 15, 35, 85, 180 $ with $ N=3 $ and $ m=5, 15, 35 $ with $ N=5 $.}
	\label{tabel-5}
\end{table}
For different values of $ \alpha $, $ m $ and $ N $, in above the corresponding radii  are represented on $ x $-axis cut by the increasing curves all have asymptotes at $ x=1 $ have been shown in Figures 3-6.
\begin{thm}\label{th-2.28}
	Let $ f\in \mathcal{P}^{0}_{\mathcal{H}}(\alpha) $ be given by \eqref{e-1.11a}. Then 
	\begin{align}\label{e-2.29}
		r+|h(r)|^p+\sum_{n=2}^{\infty}(|a_n|+|b_n|)r^{n}\leq d\left(f(0),\partial f(\mathbb{D})\right),\; \text{for}\;\;  r\leq r_{p}(\alpha),
	\end{align}
	where $ r_{p}(\alpha) $ is the smallest root of the equation
	\begin{equation}\label{e-2.30}
		r^p+r-1-2(1-\alpha)\left(r-1+\ln\,(2-2r)\right)=0 \;\; \mbox{in}\;\; (0,1).
	\end{equation}
	The radius $ r_{p}(\alpha) $ is the best possible.
\end{thm}
\vspace{1in}
\begin{figure}[!htb]
	\begin{center}
		\includegraphics[width=0.40\linewidth]{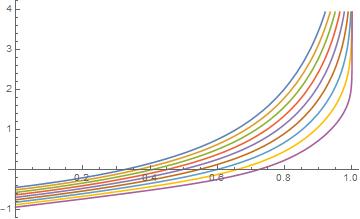}\;\;\;\;\;\includegraphics[width=0.40\linewidth]{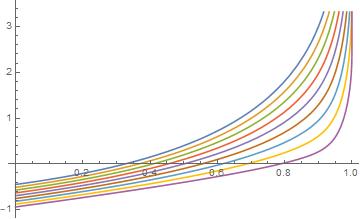}
	\end{center}
	\caption{The roots $ r_{_{_7}}(\alpha) $ and $ r_{_{35}}(\alpha) $ of the equation $ 	r^p+r-1-2(1-\alpha)\left(r-1+\ln\,(2-2r)\right)=0 $ in $ (0,1) $.}
\end{figure}
\begin{table}[ht]
	\centering
	\begin{tabular}{|l|l|l|l|l|l|l|l|l|l|}
		\hline
		$\alpha$& $0.1$&$0.2$& $0.3$& $ 0.4 $& $0.5$&$0.6$& $0.7$& $ 0.8 $& $ 0.9 $ \\
		\hline
		$r_7(\alpha)$& $0.3249 $&$0.3653 $& $0.4069$& $0.4503$& $0.4963 $&$0.5456 $& $0.5992$& $0.6579$ &$ 0.7231 $\\
		\hline
		$r_{35}(\alpha)$& $0.3251 $&$0.3657 $& $0.4078$& $0.4522$& $0.5000 $&$0.5529 $& $0.6136$& $0.6872$ &$ 0.7867 $\\
		\hline
	\end{tabular}
	\vspace{1mm}
	\caption{The roots $r_{_7}(\alpha)$ and $r_{_{35}}(\alpha)$ of the equation $ 	r^p+r-1-2(1-\alpha)\left(r-1+\ln\,(2-2r)\right)=0 $ for $ \alpha\in[0,1) $.}
	\label{tabel-6}
\end{table}

\section{Improved Bohr Radius for the class $ \mathcal{P}^{0}_{\mathcal{H}}(\alpha) $ }
In $ 2020 $, Kayumov and Ponnusamy \cite{kayumov-2018-c} have obtained several improved versions of Bohr inequality for analytic functions.
\begin{thm} \cite{kayumov-2018-c} \label{him-thm-3.1}
	Let $f(z)=\sum_{n=0}^{\infty} a_{n}z^{n}$ be analytic in $\mathbb{D}$, $|f(z)|\leq 1$ and $S_{r}$ denotes the image of the subdisk $|z|<r$ under mapping $f$. Then 
	\begin{equation}
		B_{1}(r):=\sum_{n=0}^{\infty}|a_n|r^n+ \frac{16}{9} \left(\frac{S_{r}}{\pi}\right) \leq 1 \quad \mbox{for} \quad r \leq \frac{1}{3}
	\end{equation}
	and the numbers $1/3$, $16/9$ cannot be improved. Moreover, 
	\begin{equation}
		B_{2}(r):=|a_{0}|^{2}+\sum_{n=1}^{\infty}|a_n|r^n+ \frac{9}{8} \left(\frac{S_{r}}{\pi}\right) \leq 1 \quad \mbox{for} \quad r \leq \frac{1}{2}
	\end{equation}
	and the numbers $1/2$ and $9/8$ cannot be improved.	
\end{thm}

\begin{thm}\cite{kayumov-2018-c} \label{him-thm-3.4}
	Let $f(z)=\sum_{n=0}^{\infty} a_{n}z^{n}$ be analytic in $\mathbb{D}$ and $|f(z)|\leq 1$. Then 
	\begin{equation}
		|a_{0}|+ \sum_{n=1}^{\infty}\bigg(|a_n|+\frac{1}{2}|a_n|^2\bigg)r^n\leq 1 \quad \mbox{for} \quad r \leq \frac{1}{3}
	\end{equation}
	and the numbers $1/3$ and $1/2$ cannot be improved.
\end{thm}
The primary object of this section is to generalize the harmonic versions of Theorem \ref{him-thm-3.1} and Theorem \ref{him-thm-3.4} for the class $\mathcal{P}^{0}_{\mathcal{H}}(\alpha) $. It is interesting to investigate Theorem \ref{him-thm-3.1}, when $ {S_{r}}/{\pi} $ has certain power. Therefore in order to generalize Theorem \ref{him-thm-3.1}, we consider a $ { N^{\rm th}} $ degree polynomial in $ S_r/\pi $ of the form 
\begin{equation*}
	P\left(\frac{S_{r}}{\pi}\right)=\left( \frac{S_{r}}{\pi}\right)^N+\left(\frac{S_{r}}{\pi}\right)^{N-1}+\cdots+\frac{S_{r}}{\pi}. 
\end{equation*} 

\begin{thm}\label{th-2.13}
	Let $ f\in \mathcal{P}^{0}_{\mathcal{H}}(\alpha) $ be given by \eqref{e-1.11a}. Then 
	\begin{equation}\label{e-2.14}
		r+\sum_{n=2}^{\infty}(|a_n|+|b_n|)r^n+P\left(\frac{S_{r}}{\pi}\right)\leq d\left(f(0),\partial f(\mathbb{D})\right)
	\end{equation} 
	for $ r\leq r_{_N}(\alpha) $, where $ P(w)=w^N+w^{N-1}+\cdots+w $, a polynomial in $ w $ of degree $ N-1 $, and $ r_{_N}(\alpha)\in (0,1) $ is the smallest root of the equation
	\begin{equation}\label{e-2.15}
		r-1-2(1-\alpha)\left(r-1+\ln\,(2-2r)\right)+P\left(r^2-4(1-\alpha)^2(r^2+\ln\;(1-r^2))\right)=0.
	\end{equation}
	Th constant $ r_{_N}(\alpha) $ is the best possible.
\end{thm}

\begin{table}[ht]
	\centering
	\begin{tabular}{|l|l|l|l|l|l|l|l|l|l|}
		\hline
		$\alpha$& $0.1$&$0.2$& $0.3$& $ 0.4 $& $0.5$&$0.6$& $0.7$& $ 0.8 $& $ 0.9 $ \\
		\hline
		$r_2(\alpha)$& $0.2734 $&$0.3027 $& $0.3320$& $0.3618$& $0.3923 $&$0.4241 $& $0.4574$& $0.4927$ &$ 0.5303 $\\
		\hline
		$r_{3}(\alpha)$& $0.2732 $&$0.3023 $& $0.3314$& $0.3607$& $0.3907 $&$0.4217 $& $0.4540$& $0.4878$ &$ 0.5230 $\\
		\hline
		$r_4(\alpha)$& $0.2732 $&$0.3023 $& $0.3313$& $0.3606$& $0.3905 $&$0.4213 $& $0.4533$& $0.4867$ &$ 0.5212 $\\
		\hline
		$r_5(\alpha)$& $0.2732 $&$0.3023 $& $0.3313$& $0.3606$& $0.3904 $&$0.4213 $& $0.4532$& $0.4864$ &$ 0.5208 $\\
		\hline
	\end{tabular}
	\vspace{1mm}
	\caption{The roots $ r_{_N}(\alpha) $ of equation \eqref{e-2.15} when $ N=2, 3, 4, 5 $ and $ \alpha\in[0,1) $.}
	\label{tabel-7}
\end{table}

\noindent As a consequence of Theorem \ref{th-2.13}, we obtain the following interesting corollary.
\begin{cor}\label{cor-3.4}
		Let $ f\in \mathcal{P}^{0}_{\mathcal{H}}(\alpha) $ be given by \eqref{e-1.11a}. Then 
		\begin{equation}\label{e-2.11}
			r+\sum_{n=2}^{\infty}(|a_n|+|b_n|)r^n+\left(\frac{S_{r}}{\pi}\right)\leq d\left(f(0),\partial f(\mathbb{D})\right)
		\end{equation} 
		for $ r\leq r_{\alpha} $, where $ r_{\alpha}\in (0,1) $ is the smallest root of the equation
		\begin{equation}\label{e-2.12}
			r^2+r-1-2(1-\alpha)(3-2\alpha)\left(r+\ln\,(1-r)\right)-2(1-\alpha)(\ln 2-1)=0.
		\end{equation}
		The radius $ r_{\alpha} $ is the best possible.
\end{cor}
	\begin{figure}[!htb]
	\begin{center}
		\includegraphics[width=0.55\linewidth]{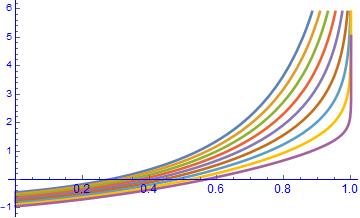}
	\end{center}
	\caption{The figure demonstrates the roots $ r(\alpha) $ when $ \alpha\in[0,1) $.}
\end{figure}
\begin{table}[ht]
	\centering
	\begin{tabular}{|l|l|l|l|l|l|l|l|l|l|}
		\hline
		$\alpha$& $0.1$&$0.2$& $0.3$& $ 0.4 $& $0.5$&$0.6$& $0.7$& $ 0.8 $& $ 0.9 $ \\
		\hline
		$r(\alpha)$& $0.2322 $&$0.2635 $& $0.2967$& $0.3323$& $0.3707 $&$0.4125 $& $0.4579$& $0.5074$ &$ 0.5610 $\\
		\hline
	\end{tabular}
	\vspace{1mm}
	\caption{The roots $ r(\alpha) $ of the equation \eqref{e-2.12} when $ \alpha\in[0,1) $.}
	\label{tabel-8}
\end{table}
 The polylogarithm function is defined by a power series in $ z $, which is also a Dirichlet series in $ s $. That is
\begin{equation*}
	Li_s(z)=\sum_{n=1}^{\infty}\frac{z^n}{n^s}=z+\frac{z^2}{2^s}+\frac{z^3}{3^s}+\cdots,
\end{equation*} 
is valid for arbitrary complex order $ s $ and for all complex arguments $ z $ with $ |z|<1. $ Therefore, the dilogarithm function is denoted by $ Li_2(z) $, is a particular case of the polylogarithm. The following theorem is the generalization of the harmonic version of Theorem \ref{him-thm-3.4} by considering the right hand side $ d\left(f(0),\partial f(\mathbb{D})\right) $ instead of $ 1 $.
\begin{thm}\label{th-2.16}
	Let $ f\in \mathcal{P}^{0}_{\mathcal{H}}(\alpha) $  be given by \eqref{e-1.11a}. Then 
	\begin{equation}\label{e-2.17}
		r+\sum_{n=2}^{\infty}\left(|a_n|+|b_n|+(|a_n|+|b_n|)^2\right)r^n\leq d\left(f(0),\partial f(\mathbb{D})\right),\; \text{for}\;\;  r\leq r_{\alpha},
	\end{equation} 
	where $ r_{\alpha}\in (0,1) $ is the smallest root of the equation
	\begin{equation}\label{e-2.18}
		r-1-2(1-\alpha)\left(r-1+\ln\,(2-2r)\right)+4(1-\alpha^2)(Li_{2}(r)-r)=0,
	\end{equation} where $ Li_{2}(z) $ is a dilogarithm. The constant $ r_{\alpha} $ is best possible.
\end{thm}

\section{Refined Bohr Radius for the class $ \mathcal{P}^{0}_{\mathcal{H}}(\alpha) $ }
In $ 2020 $, Ponnusamy {\it et al.} \cite{Ponnusamy-Vijaya-ResultsMath-2020} established the following refined Bohr inequality by applying a refined version of the coefficient inequalities.
\begin{thm} \cite{Ponnusamy-Vijaya-ResultsMath-2020} \label{him-thm-4.1}
	Let $f(z)=\sum_{n=0}^{\infty} a_{n}z^{n}$ be analytic in $\mathbb{D}$ and $|f(z)|\leq 1$. Then
	$$
	\sum_{n=0}^{\infty} |a_{n}|r^{n}+ \left(\frac{1}{1+|a_{0}|}+\frac{r}{1-r}\right) \sum_{n=1}^{\infty}|a_n|^2r^{2n}\leq 1
	$$
	for $r\leq 1/(2+|a_{0}|)$ and the numbers $1/(1+|a_{0}|)$ and $1/(2+|a_{0}|)$ cannot be improved. Moreover, 
	$$
	|a_{0}|^{2}+ \sum_{n=1}^{\infty} |a_{n}|r^{n}+ \left(\frac{1}{1+|a_{0}|}+\frac{r}{1-r}\right) \sum_{n=1}^{\infty}|a_n|^2r^{2n}\leq 1
	$$
	for $r\leq 1/2$ and the numbers $1/(1+|a_{0}|)$ and $1/2$ cannot be improved.
\end{thm}
\noindent Next we prove the harmonic analogue of Theorem \ref{him-thm-4.1}.
\begin{thm}\label{th-2.19}
	Let $ f\in \mathcal{P}^{0}_{\mathcal{H}}(\alpha) $ be given by \eqref{e-1.11a}. Then 
	\begin{align}
		\label{e-2.20}
		& r+\sum_{n=2}^{\infty}(|a_n|+|b_n|)r^n+\frac{1}{1-r^N}\sum_{n=2}^{\infty}n(|a_n|+|b_n|)^2r^{2n}\\&\nonumber\leq d\left(f(0),\partial f(\mathbb{D})\right)\; \text{for}\;\;  r\leq r_{_N}(\alpha),
	\end{align}
	where $ r_{_N}(\alpha)\in (0,1) $ is the smallest root of the equation
	\begin{equation}\label{e-2.21}
		r-1-2(1-\alpha)\left(r-1+\ln\,(2-2r)+\frac{2(1-\alpha)}{1-r^N}(r^2+\ln\,(1-r^2))\right)=0.
	\end{equation} Here $ r_{_N}(\alpha) $ is the best possible.
\end{thm}

\begin{figure}[!htb]
	\begin{center}
		\includegraphics[width=0.48\linewidth]{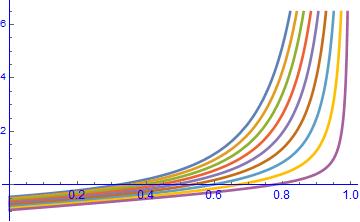}\;\;\;\;\;\includegraphics[width=0.48\linewidth]{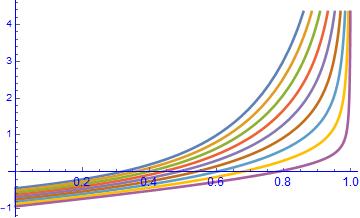}
	\end{center}
	\caption{The roots $ r_{2}(\alpha) $ and $ r_{25}(\alpha) $ when $ \alpha\in [0,1) $ have been shown in Figure 9.}
\end{figure}

\begin{table}[ht]
	\centering
	\begin{tabular}{|l|l|l|l|l|l|l|l|l|l|}
		\hline
		$\alpha$& $0.1$&$0.2$& $0.3$& $ 0.4 $& $0.5$&$0.6$& $0.7$& $ 0.8 $& $ 0.9 $ \\
		\hline
		$r_2(\alpha)$& $0.3148 $&$0.3527 $& $0.3920$& $0.4338$& $0.4793 $&$0.5304 $& $0.5904$& $0.6651$ &$ 0.7693 $\\
		\hline
		$r_{25}(\alpha)$& $0.3158 $&$0.3542 $& $0.3942$& $0.4368$& $0.4835 $&$0.5361 $& $0.5977$& $0.6741$ &$ 0.7792 $\\
		\hline
	\end{tabular}
	\vspace{1mm}
	\caption{In this table,  we obtained the roots of $ r_{_N}(\alpha) $ for $ N=2 $ and $ 25 $ when $ \alpha\in[0,1) $.}
	\label{tabel-3}
\end{table}
In the following, we prove two interesting results which are harmonic analogue of refined Bohr inequality.
\begin{thm}\label{th-2.22}
	Let $ f\in \mathcal{P}^{0}_{\mathcal{H}}(\alpha) $  be given by \eqref{e-1.11a}. Then 
	\begin{align}
		\label{e-2.23}
		& r+\sum_{n=2}^{\infty}(|a_n|+|b_n|)r^n+\left(\frac{1}{1+|a_2|+|b_2|}+\frac{r^m}{1-r^m}\right)\sum_{n=3}^{\infty}n^{m-1}(|a_n|+|b_n|)^mr^{mn}\\&\nonumber\leq d\left(f(0),\partial f(\mathbb{D})\right),\; \text{for}\;\;  r\leq r_{{m}}(\alpha),
	\end{align}
	where $ r_{{m}}(\alpha)\in (0,1) $ is the smallest root of the equation
\begin{align}
	\label{e-2.24} &
	r-1-2(1-\alpha)(r-1+\ln\,(2-2r))\\&\nonumber-2^m(1-\alpha)^m\left(\frac{1}{1+|a_2|+|b_2|}+\frac{r^m}{1-r^m}\right)\left(r^m+\frac{r^{2m}}{2}+\ln\,(1-r^m)\right)=0.
\end{align}
The constant $ r_{{m}}(\alpha) $ is the best possible.
\end{thm}
\begin{thm}\label{th-2.25}
	Let $ f\in \mathcal{P}^{0}_{\mathcal{H}}(\alpha) $ be given by \eqref{e-1.11a}. Then 
	\begin{align}\label{e-2.26}
		& r+\frac{\left(1-(1+|a_2|+|b_2|-(|a_2|+|b_2|)^2)\right)r}{1-(|a_2|+|b_2|)r}+\sum_{n=3}^{\infty}(|a_n|+|b_n|)r^{n}\\&\nonumber\leq d\left(f(0),\partial f(\mathbb{D})\right)\; \text{for}\;\;  r\leq r_{\alpha},
	\end{align}
	where $ r_{\alpha}\in (0,1) $ is the smallest root of the equation
\begin{equation}\label{e-2.27}
r-\frac{\left(1-(|a_2|+|b_2|)^2\right)r}{1-(|a_2|+|b_2|)r}-2(1-\alpha)\left(r+\frac{r^2}{2}-1+\ln\,(2-2r)\right)=0.
\end{equation}
The constant $ r_{\alpha} $ is the best possible.
\end{thm}
\section{Bohr-Type Inequality for the class $ \mathcal{P}^{0}_{\mathcal{H}}(\alpha) $}
\noindent We now prove the following Bohr-type inequality for the class of functions $ \mathcal{P}^{0}_{\mathcal{H}}(\alpha) $.
\begin{thm}\label{th-2.31}
	Let $ f\in \mathcal{P}^{0}_{\mathcal{H}}(\alpha) $  be given by \eqref{e-1.11a}. Then 
	\begin{align}\label{e-2.32}
		|f(z)|+\sqrt{|{J}_f(z)|}r+\sum_{n=N}^{\infty}(|a_n|+|b_n|)r^{n}\leq d\left(f(0),\partial f(\mathbb{D})\right)\; \text{for}\;\;  r\leq r_{_N}(\alpha),
	\end{align}
	where $ r_{_N}(\alpha)\in (0,1) $ is the smallest root of the equation

	\begin{align}	\label{e-2.33} &
		r-1-2(1-\alpha)\left(2r-1+\frac{r^2}{2}+\cdots+\frac{r^{N-1}}{N-1}+\ln 2+2\ln(1-r)\right)\\&\quad\quad+\left(\alpha+(1-\alpha)\left(\frac{1+r}{1-r}\right)\right)r=0\nonumber.
	\end{align}
	The radius $ r_{_N}(\alpha) $ is the best possible.
\end{thm}
%______________________________________________________________________%

\section{Proof of the main results}	
\begin{pf}[\bf Proof of Theorem \ref{th-2.1}]
	Let $f \in \mathcal{P}^{0}_{\mathcal{H}}(\alpha)$ be given by \eqref{e-1.11a}. Then from Theorem \ref{th-2.1}, we have
	
	\begin{equation} \label{e-3.1}
		|f(z)|\geq |z|+(1-\alpha)\sum\limits_{n=2}^{\infty}  \dfrac{2(-1)^{n-1}}{n} |z|^{n} \quad \mbox{for } \quad |z|<1.
	\end{equation}
The Euclidean distance between $f(0)$ and the boundary of $f(\mathbb{D})$ is given by 
\begin{equation} \label{e-3.2}
d(f(0), \partial f(\mathbb{D}))= \liminf \limits_{|z|\rightarrow 1} |f(z)-f(0)|.
\end{equation}
	Since $f(0)=0$, from \eqref{e-1.10} and \eqref{e-3.2} we obtain 
	\begin{equation} \label{e-3.3}
		d(f(0), \partial f(\mathbb{D})) \geq 1+\sum\limits_{n=2}^{\infty}  2(1-\alpha) \dfrac{(-1)^{n-1}}{n}.
	\end{equation}
Using Lemmas \ref{lem-1.8} and \ref{lem-1.9}, for $ |z|=r_{_N}(\alpha) $, we obtain 
\begin{align}\label{e-3.4} &
	|f(z)|+\sum_{n=N}^{\infty}(|a_n|+|b_n|)r^n \\&\nonumber\leq r+\sum_{n=2}^{\infty}\frac{2(1-\alpha)r^n}{n}+\sum_{n=N}^{\infty}\frac{2(1-\alpha)r^n}{n}\\&=\nonumber r-2(1-\alpha)(r+\ln(1-r))-2(1-\alpha)\left(\ln(1-r)+\sum_{n=1}^{N-1}\frac{r^n}{n}\right)\\&=\nonumber r-2(1-\alpha)\left(r+\ln((1-r)^2)+\sum_{n=1}^{N-1}\frac{r^n}{n}\right).
\end{align}
It is easy to see that 
\begin{equation}\label{e-3.5}
	r-2(1-\alpha)\left(r+\ln((1-r)^2)+\sum_{n=1}^{N-1}\frac{r^n}{n}\right)\leq 1+2(1-\alpha)(\ln 2-1)
\end{equation}
for $ r\leq r_{_N}(\alpha) $, where $ r_{_N}(\alpha) $ is the smallest root of
\begin{equation*}
	r-1-2(1-\alpha)\left(r-1+\ln(2(1-r)^2)+\sum_{n=1}^{N-1}\frac{r^n}{n}\right)= 0
\end{equation*}
in $ (0,1) $.
Let $ H_1 : [0,1)\rightarrow \mathbb{R} $ be defined by 
$$ H_1(r):=r-1-2(1-\alpha)\left(r-1+\ln(2(1-r)^2)+\sum_{n=1}^{N-1}\frac{r^n}{n}\right).
$$
 The existence of a root $ r_N(\alpha) $ is ensured by the following fact that $ H_1 $ is a continuous function with the properties $ H_1(0)=-1-2(1-\alpha)(\ln 2-1)<0 $ and $ \displaystyle\lim_{r\rightarrow 1}H_1(r)=+\infty. $ Let $ r_N(\alpha) $ to be the smallest root of $ H_1(r)=0 $ in $ (0,1) $. Therefore, we have $ H_1(r_N(\alpha))=0 $. That is 
 \begin{align}\label{e-3.6} &
 		r_N(\alpha)-1-2(1-\alpha)\left(r_N(\alpha)-1+\ln(2(1-r_N(\alpha))^2)+\sum_{n=1}^{N-1}\frac{r^n_N(\alpha)}{n}\right)=0.
 \end{align}
 
  In view of \eqref{e-3.3}, \eqref{e-3.4} and \eqref{e-3.5} for $ |z|=r\leq r_N(\alpha) $, it follows that
 $$ 
 |f(z)|+\sum_{n=N}^{\infty}(|a_n|+|b_n|)r^n\leq d(f(0),\partial f(\mathbb{D})).
  $$ 
  In order to show that the constant $ r_{_N}(\alpha) $ is the best possible constant, we consider the following function $ f=f_{\alpha} $ by 
  \begin{equation}\label{e-6.7}
  	f_{\alpha}(z)=z+\sum_{n=2}^{\infty}\frac{2(1-\alpha)z^n}{n}.
  \end{equation}
  It is easy to show that $ f_{\alpha}\in\mathcal{P}^{0}_{\mathcal{H}}(\alpha) $. For $ f=f_{\alpha} $, it can be seen that 
  \begin{equation}\label{e-3.7}
d(f(0),\partial f(\mathbb{D}))=1+2(1-\alpha)(\ln 2 -1).
  \end{equation}
For $ f=f_{\alpha} $ and $ |z|=r_N(\alpha) $, a simple computation using \eqref{e-3.6} and \eqref{e-3.7} shows that 
\begin{align*} &
		|f(z)|+\sum_{n=N}^{\infty}(|a_n|+|b_n|)r^n\\ &= r_{_N}(\alpha)+\sum_{n=2}^{\infty}\frac{2(1-\alpha)(r_{_N}(\alpha))^n}{n}+\sum_{n=N}^{\infty}\frac{2(1-\alpha)(r_{_N}(\alpha))^n}{n}\\&=
		r_{_N}(\alpha)-2(1-\alpha)\left(r_{_N}(\alpha)+\ln((1-r_{_N}(\alpha))^2)+\sum_{n=1}^{N-1}\frac{r^n_{_N}(\alpha)}{n}\right)\\&\nonumber= 1+2(1-\alpha)(\ln 2-1)\\&=d(f(0),\partial f(\mathbb{D})).
\end{align*}
Therefore, the radius $ r_{_N}(\alpha) $ is the best possible. This completes the proof.
\end{pf}

\begin{pf}[\bf Proof of Theorem \ref{th-2.4}]
Let $f \in \mathcal{P}^{0}_{\mathcal{H}}(\alpha)$ be given by \eqref{e-1.11a}.  Then in view of Lemmas \ref{lem-1.8} and \ref{lem-1.9}, for $ |z|=r $, we obtain 
\begin{align}\label{e-3.10}
		|f(z)|^2+\sum_{n=N}^{\infty}|a_n||z|^n &\leq \left(r+\sum_{n=2}^{\infty}\frac{2(1-\alpha)r^n}{n}\right)^2+\sum_{n=N}^{\infty}\frac{2(1-\alpha)r^n}{n}.
\end{align}
A simple computation shows that
\begin{equation}\label{e-3.11}
	\left(r+\sum_{n=2}^{\infty}\frac{2(1-\alpha)r^n}{n}\right)^2+\sum_{n=N}^{\infty}\frac{2(1-\alpha)r^n}{n}\leq 1+2(1-\alpha)(\ln 2-1)
\end{equation}
for $ r\leq r_N(\alpha) $, where $ r_N(\alpha) $ is the smallest root of $ H_2(r)=0 $ in $ (0,1) $, where $ H_2 : [0,1)\rightarrow \mathbb{R} $ is defined by 
$$ H_2(r)=\bigg(r-2(1-\alpha)(r+\ln(1-r))\bigg)^2-2(1-\alpha)\left(\ln (2-2r)-1+\sum_{n=1}^{N-1}\frac{r^n}{n}\right)- 1.
$$
Then $ H_2 $ is a continuous function with $ H_2(0)=-1-2(1-\alpha)(\ln 2-1)<0 $ and $ \displaystyle\lim_{r\rightarrow 1}H_2(r)=+\infty. $  Therefore, $ H_2(r)=0 $ has a root in $ (0,1) $ and we choose the smallest root to be $ r_{_N}(\alpha) $. Therefore, we have $ H_2(r_{_N}(\alpha))=0 $. That is
\begin{align}\label{e-3.12} &
	\bigg(r_N(\alpha)-2(1-\alpha)(r_N(\alpha)+\ln(1-r_N(\alpha)))\bigg)^2\\&\nonumber -2(1-\alpha)\left(\ln (2-2r_N(\alpha))-1+\sum_{n=1}^{N-1}\frac{r^n_N(\alpha)}{n}\right)- 1=0.
\end{align}

\noindent Using \eqref{e-3.3}, \eqref{e-3.10} and \eqref{e-3.11} for $ |z|=r\leq r_{_N}(\alpha) $, we obtain
$$ 
|f(z)|^2+\sum_{n=N}^{\infty}|a_n|r^n\leq d(f(0),\partial f(\mathbb{D})).
$$ 
In order to show that $ r_{_N}(\alpha) $ is the best possible, we consider the function $ f=f_{\alpha} $ defined by \eqref{e-6.7}. For $ f=f_{\alpha} $ and $ |z|=r_{_N}(\alpha) $, a simple computation using \eqref{e-3.7} and \eqref{e-3.12} shows that 
\begin{align*} 
	|f(z)|^2+\sum_{n=N}^{\infty}|a_n|r^n&= \left(r_{_N}(\alpha)+\sum_{n=2}^{\infty}\frac{2(1-\alpha)(r_{_N}(\alpha))^n}{n}\right)^2+\sum_{n=N}^{\infty}\frac{2(1-\alpha)(r_{_N}(\alpha))^n}{n}\\&\nonumber= 1+2(1-\alpha)(\ln 2-1)\\&=d(f(0),\partial f(\mathbb{D})).
\end{align*}
Therefore, the radius $ r_{_N}(\alpha) $ is the best possible.	
\end{pf}

\begin{pf}[\bf Proof of Theorem \ref{th-2.7}]
	Let $f \in \mathcal{P}^{0}_{\mathcal{H}}(\alpha)$ be given by \eqref{e-1.11a}. Using Lemmas \ref{lem-1.8} and \ref{lem-1.9} for $ |z|=r $, we obtain 
	\begin{align}\label{e-6.12}
		|f(z^m)|+\sum_{n=N}^{\infty}|a_n||z|^n &\leq r^m+\sum_{n=2}^{\infty}\frac{2(1-\alpha)r^{mn}}{n}+\sum_{n=N}^{\infty}\frac{2(1-\alpha)r^n}{n}.
	\end{align}
A simple computation shows that
\begin{align}
	\label{e-6.13} &
r^m+\sum_{n=2}^{\infty}\frac{2(1-\alpha)(r^m)^{n}}{n}+\sum_{n=N}^{\infty}\frac{2(1-\alpha)r^n}{n}\\&=\nonumber r^m-2(1-\alpha)(r^m+\ln(1-r^m))-2(1-\alpha)\left(r+\ln (1-r)+\sum_{n=2}^{N-1}\frac{r^n}{n}\right)\\&=r^m-2(1-\alpha)\left(r^m+r+\ln (1-r)(1-r^m)+\sum_{n=2}^{N-1}\frac{r^n}{n}\right)\nonumber\\& \nonumber\leq 1+2(1-\alpha)(\ln 2-1)
\end{align}
	for $ r\leq r{_{m,N}}(\alpha) $, where $ r_{_{m,N}}(\alpha) $ is the smallest root of $ H_3(r)=0 $ in $ (0,1) $, where $ H_3 : [0,1)\rightarrow \mathbb{R} $ is defined by 
	$$ H_3(r):=r^m-1-2(1-\alpha)\left(r^m+r-1+\ln (2-2r)(1-r^m)+\sum_{n=2}^{N-1}\frac{r^n}{n}\right).
	$$
In view of the same line of argument as in the proof of Theorem \ref{th-2.1}, we can show that $ H_3(r)=0 $ has a root in $ (0,1) $ and we choose $ r_{_{m,N}}(\alpha) $ to be the smallest root of $ H_3(r) $. Therefore, we have $ H_3(r_{_{m,N}}(\alpha))=0 $. That is
	\begin{align}\label{e-6.14} 
		r_{_{m,N}}(\alpha)-1-2(1-\alpha)G_{m,N}(r)=0,
	\end{align}
    where 
    \begin{equation*}
    	G_{m,N}=\left(r_{_{m,N}}(\alpha)+r_{_N}(\alpha)-1+\ln\left((2-2r_N(\alpha))(1-r_{_{m,N}}(\alpha))\right) +\sum_{n=2}^{N-1}\frac{r^{n}_{_{m,N}}(\alpha)}{n}\right).
    \end{equation*}
	In view of \eqref{e-3.3}, \eqref{e-6.12} and \eqref{e-6.13} for $ |z|=r\leq r_{_{m,N}}(\alpha) $, we obtain
	$$ 
	|f(z^m)|+\sum_{n=N}^{\infty}|a_n|r^n\leq d(f(0),\partial f(\mathbb{D})).
	$$ 
	To show that the radius $ r_{_{m,N}}(\alpha) $ is the best possible, we consider the function $ f=f(\alpha) $ defined by \eqref{e-6.7}. 	For $ f=f_{\alpha} $ and $ |z|=r_{_{m,N}}(\alpha) $, a simple calculation using \eqref{e-3.7} and \eqref{e-6.14} shows that 
	\begin{align*} &
|f(z^m)|+\sum_{n=N}^{\infty}|a_n|r^n\\&= r^m_{_{m,N}}(\alpha)-2(1-\alpha)\left(r^m_{_{m,N}}(\alpha)+r+\ln \bigg((1-r_{_{m,N}}(\alpha))(1-r^m_{_{m,N}}(\alpha))\bigg)+\sum_{n=2}^{N-1}\frac{r^n_{_{m,N}}(\alpha)}{n}\right)\\&\nonumber= 1+2(1-\alpha)(\ln 2-1)=d(f(0),\partial f(\mathbb{D})).
	\end{align*}
	Hence, the radius $ r_{_{m,N}}(\alpha) $ is the best possible. This completes the proof.
\end{pf}

\begin{pf}[\bf Proof of Theorem \ref{th-2.28}]
	Let $f \in \mathcal{P}^{0}_{\mathcal{H}}(\alpha)$ be given by \eqref{e-1.11a}. Applying Lemmas \ref{lem-1.8} and \ref{lem-1.9} for $ |z|=r $, we obtain 
	\begin{align}\label{e-6.15}
		r+|h(r)|^p+\sum_{n=2}^{\infty}(|a_n|+|b_n|)r^{n} &\leq r+|h(r)|^p+\sum_{n=2}^{\infty}\frac{2(1-\alpha)r^{n}}{n}.
	\end{align}
	It is not difficult to show that
	\begin{align}
		\label{e-6.16} 
		r+|h(r)|^p+\sum_{n=2}^{\infty}\frac{2(1-\alpha)r^{n}}{n}&=r^p+r-2(1-\alpha)(r+\ln(1-r))\\& \nonumber\leq 1+2(1-\alpha)(\ln 2-1)
	\end{align}
	for $ r\leq r_p(\alpha) $, where $ r_p(\alpha) $ is the smallest root of $ H_4(r)=0 $	in $ (0,1) $ and $ H_{4} : [0,1)\rightarrow \mathbb{R} $ is defined by 
	$$ H_{4}(r):=r^p+r-1-2(1-\alpha)(r-1+\ln(2-2r)).
	$$
	By the same argument used in the proof of Theorem \ref{th-2.1}, we can show that $ H_4(r) $ has a root in $ (0,1) $ and we choose $ r_p(\alpha) $ to be the smallest root of $ H_4(r) $. Therefore, $ H_{4}(r_p(\alpha))=0 $. That is  
	\begin{align}\label{e-6.17} 
		r^p_p(\alpha)+r_p(\alpha)-1-2(1-\alpha)(r_p(\alpha)-1+\ln(2-2r_p(\alpha)))=0,
	\end{align}
	Using \eqref{e-3.3}, \eqref{e-6.15} and \eqref{e-6.16} for $ |z|=r\leq r_p(\alpha) $, we obtain
	$$ 
	r+|h(r)|^p+\sum_{n=2}^{\infty}(|a_n|+|b_n|)r^{n}\leq d(f(0),\partial f(\mathbb{D})).
	$$ 
	In order to show that $ r_p(\alpha) $ is the best possible, we consider the function $ f=f_{\alpha} $ defined by \eqref{e-6.7}.	For $ f=f_{\alpha} $ and $ |z|=r_p(\alpha) $, a simple calculation using \eqref{e-3.7} and \eqref{e-6.17} shows that  \vspace{-0.4in}
	\begin{align*}&
		r_p(\alpha)+|h(r_p(\alpha))|^p+\sum_{n=2}^{\infty}(|a_n|+|b_n|)r^{n}\\&= r^p_p(\alpha)+r_p(\alpha)-2(1-\alpha)(r_p(\alpha)+\ln(1-r_p(\alpha)))\\&\nonumber= 1+2(1-\alpha)(\ln 2-1)\\&=d(f(0),\partial f(\mathbb{D})).
	\end{align*}
	Therefore, the radius $ r_p(\alpha) $ is the best possible.	This completes the proof.	
\end{pf}

 \begin{pf}[\bf Proof of Theorem \ref{th-2.13}]
	Let $f \in \mathcal{P}^{0}_{\mathcal{H}}(\alpha)$ be given by \eqref{e-1.11a}.
For the analytic functions $ h $ and $ g $, the area of the disk $ |z|<r $ under the harmonic map $ f $ is $ S_r $ is given by
\begin{align}\label{e-6.18a}
	S_r=\iint\limits_{D_r}&\left(|h^{\prime}(z)|^2-|g^{\prime}(z)|^2\right)dxdy,\\\label{e-6.18b} \frac{1}{\pi}\iint\limits_{D_r}&|h^{\prime}(z)|^2dxdy=\sum_{n=1}^{\infty}n|a_n|^2r^{2n},\\\label{e-6.18c} \frac{1}{\pi}\iint\limits_{D_r}&|g^{\prime}(z)|^2dxdy=\sum_{n=2}^{\infty}n|b_n|^2r^{2n}.
\end{align}
\noindent Therefore, in view of \eqref{e-6.18a}, \eqref{e-6.18b} and \eqref{e-6.18c} and Lemma \ref{lem-1.8}, we obtain
\begin{align*}
	\frac{S_r}{\pi}&=\frac{1}{\pi}\iint\limits_{D_r}\left(|h^{\prime}(z)|^2-|g^{\prime}(z)|^2\right)dxdy\\&=r^2+\sum_{n=2}^{\infty}|a_n|^2r^{2n}-\sum_{n=2}^{\infty}n|b_n|^2r^{2n}\\&= r^2+\sum_{n=2}^{\infty}n\left(|a_n|+|b_n|\right)\left(|a_n|-|b_n|\right)r^{2n}\\&\leq r^2+\sum_{n=2}^{\infty}\frac{4(1-\alpha)^2r^{2n}}{n^2}\\&=r^2-4(1-\alpha)^2\left(r^2+\ln (1-r^2)\right).
\end{align*}

\noindent Using Lemmas \ref{lem-1.8} and \ref{lem-1.9} for $ |z|=r $, we obtain 
\begin{align}\label{e-6.18}
	& r+\sum_{n=2}^{\infty}(|a_n|+|b_n|)r^n+P\left(\frac{S_{r}}{\pi}\right)\\&\nonumber\leq r+\sum_{n=2}^{\infty}\frac{2(1-\alpha)r^n}{n}+P(r^2-4(1-\alpha)^2\left(r^2+\ln (1-r^2)\right))\\&=\nonumber r-2(1-\alpha)(r+\ln(1-r))+P(r^2-4(1-\alpha)^2\left(r^2+\ln (1-r^2)\right))\\&\nonumber\leq 1+2(1-\alpha)(\ln 2-1)
\end{align}
for $ r\leq r(\alpha) $, where $ r(\alpha) $ is the smallest root of $ H_5(r)=0 $ in $ (0,1) $ where
 $ H_5 : [0,1)\rightarrow \mathbb{R} $ be defined by 
$$ H_5(r):=r-1-2(1-\alpha)(r-1+\ln(2-2r))+P(r^2-4(1-\alpha)^2\left(r^2+\ln (1-r^2)\right)).
$$ Clearly, $ H_5(r(\alpha))=0 $. That is 
\begin{align}\label{e-6.19} 
	r(\alpha)-1-2(1-\alpha)(r(\alpha)-1+\ln(2-2r(\alpha)))+P(G(r,\alpha))=0,
\end{align} where
\begin{equation*}
	G(r,\alpha)=r(\alpha)^2-4(1-\alpha)^2\left(r(\alpha)^2+\ln (1-r(\alpha)^2)\right).
\end{equation*}
From \eqref{e-3.3}, \eqref{e-6.18} and \eqref{e-6.19}, $ |z|=r\leq r(\alpha) $, we obtain
$$ 
r+\sum_{n=2}^{\infty}(|a_n|+|b_n|)r^n+P\left(\frac{S_{r}}{\pi}\right)\leq d(f(0),\partial f(\mathbb{D})).
$$ 
To show that the radius $ r(\alpha) $ is the best possible, we consider the function defined by \eqref{e-6.7}. For $ f=f_{\alpha} $ and $ |z|=r(\alpha) $, a simple computation using \eqref{e-3.7} and \eqref{e-6.19} shows that 
\begin{align*} 
	r+\sum_{n=2}^{\infty}(|a_n|+|b_n|)r^n+P\left(\frac{S_{r}}{\pi}\right)&= r-2(1-\alpha)(r+\ln(1-r))+P(G_(r,\alpha))\\&\nonumber= 1+2(1-\alpha)(\ln 2-1)=d(f(0),\partial f(\mathbb{D})).
\end{align*}
Thus, the radius $ r(\alpha) $ is the best possible. This completes the proof.	
\end{pf} 

\begin{pf}[\bf Proof of Theorem \ref{th-2.16}]
		Let $f \in \mathcal{P}^{0}_{\mathcal{H}}(\alpha)$ be given by \eqref{e-1.11a}. Using Lemmas \ref{lem-1.8} and \ref{lem-1.9}, for $ |z|=r $, we obtain 
	\begin{align}\label{e-6.20}
	\;\;\;\;	r+\sum_{n=2}^{\infty}\left(|a_n|+|b_n|+(|a_n|+|b_n|)^2\right)r^n &\leq r+\sum_{n=2}^{\infty}\left(\frac{2(1-\alpha)}{n}+\frac{4(1-\alpha)^2}{n^2}\right)r^n.
	\end{align}
	An easy computation shows that
	\begin{align}
		\label{e-6.21} &
		r+\sum_{n=2}^{\infty}\left(\frac{2(1-\alpha)}{n}+\frac{4(1-\alpha)^2}{n^2}\right)r^n\\&=r-2(1-\alpha)(r+\ln (1-r))+4(1-\alpha)^2(Li_2(r)-r)\nonumber\\& \nonumber\leq 1+2(1-\alpha)(\ln 2-1)
	\end{align}
	for $ r\leq r(\alpha) $, where $ r(\alpha) $ is the smallest root of $ H_6(r)=0 $ in $ (0,1) $ and
	$ H_6 : [0,1)\rightarrow \mathbb{R} $ be defined by 
	$$ H_6(r):=	r-2(1-\alpha)(r-1+\ln (2-2r))+4(1-\alpha)^2(Li_2(r)-r)-1.
	$$ Thus, we have $ H_6(r(\alpha))=0 $. That is
	\begin{align}\label{e-6.22} 
		r(\alpha)-2(1-\alpha)(r(\alpha)-1+\ln(2-2r(\alpha)))+4(1-\alpha)^2(Li_2(r(\alpha))-r(\alpha))-1=0.
	\end{align}
	Using \eqref{e-3.3}, \eqref{e-6.20} and \eqref{e-6.21} for $ |z|=r\leq r(\alpha) $, we obtain
	$$ 
	r+\sum_{n=2}^{\infty}\left(|a_n|+|b_n|+(|a_n|+|b_n|)^2\right)r^n\leq d(f(0),\partial f(\mathbb{D})).
	$$ 
	In order to show that $ r(\alpha) $ is the best possible, we consider the function defined by \eqref{e-6.7}. For $ f=f_{\alpha} $ and $ |z|=r(\alpha) $, a simple computation using \eqref{e-3.7} and \eqref{e-6.22} shows that 
	\begin{align*} &
	r(\alpha)+\sum_{n=2}^{\infty}\left(|a_n|+|b_n|+(|a_n|+|b_n|)^2\right)r^n(\alpha)\\&= r(\alpha)-2(1-\alpha)(r(\alpha)+\ln (1-r(\alpha))+4(1-\alpha)^2(Li_2(r(\alpha))-r(\alpha))\\&\nonumber= 1+2(1-\alpha)(\ln 2-1)\\&=d(f(0),\partial f(\mathbb{D})).
	\end{align*}
	Therefore, the radius $ r(\alpha) $ is the best possible. This completes the proof.	
\end{pf} 
\begin{pf}[\bf Proof of Theorem \ref{th-2.19}]
		Let $f \in \mathcal{P}^{0}_{\mathcal{H}}(\alpha)$ be given by \eqref{e-1.11a}. In view of Lemmas \ref{lem-1.8} and \ref{lem-1.9} for $ |z|=r $, we obtain 
	\begin{align}\label{e-6.23} &
	r+\sum_{n=2}^{\infty}(|a_n|+|b_n|)r^n+\frac{1}{1-r^N}\sum_{n=2}^{\infty}n(|a_n|+|b_n|)^2r^{2n}\\\nonumber &\leq r+\sum_{n=2}^{\infty}\frac{2(1-\alpha)r^{n}}{n}+\frac{1}{1-r^N}\sum_{n=2}^{\infty}\frac{4(1-\alpha)^2r^{2n}}{n}.
	\end{align}
	A simple computation shows that
	\begin{align}
		\label{e-6.24} &
		r+\sum_{n=2}^{\infty}\frac{2(1-\alpha)r^{n}}{n}+\frac{1}{1-r^N}\sum_{n=2}^{\infty}\frac{4(1-\alpha)^2r^{2n}}{n}\\&=\nonumber r-2(1-\alpha)(r+\ln(1-r))-\frac{4(1-\alpha)^2}{1-r^N}(r^2+\ln(1-r^2))\\& \nonumber\leq 1+2(1-\alpha)(\ln 2-1)
	\end{align}
	for $ r\leq r(\alpha) $, where $ r(\alpha) $ is the smallest root of $ H_7(r)=0 $ in $ (0,1) $ and $ H_7 : [0,1)\rightarrow \mathbb{R} $ be defined by 
	$$ H_7(r):=r-2(1-\alpha)(r-1+\ln(2-2r))-\frac{4(1-\alpha)^2}{1-r^N}(r^2+\ln(1-r^2))-1.
	$$ 
	Thus, we have $ H_7(r(\alpha))=0 $. That is 
	\begin{align}\label{e-6.25} &
	r(\alpha)-2(1-\alpha)(r-1+\ln(2-2(r(\alpha)))\\&\nonumber -\frac{4(1-\alpha)^2}{1-r^N(\alpha))}\bigg(r^2(\alpha)+\ln(1-r^2(\alpha))\bigg)-1=0.
	\end{align}
	Using \eqref{e-3.3}, \eqref{e-6.23} and \eqref{e-6.24} for $ |z|=r\leq r(\alpha) $, we obtain
	$$ 
r+\sum_{n=2}^{\infty}(|a_n|+|b_n|)r^n+\frac{1}{1-r^N}\sum_{n=2}^{\infty}n(|a_n|+|b_n|)^2r^{2n}\leq d(f(0),\partial f(\mathbb{D})).
	$$ 
	In order to show that $ r(\alpha) $ is the best possible, we consider the function $ f=f_{\alpha} $ defined by \eqref{e-6.7}.
	For $ f=f_{\alpha} $ and $ |z|=r(\alpha) $, a simple calculation using \eqref{e-3.7} and \eqref{e-6.25} shows that 
	\begin{align*}  &
		r(\alpha)+\sum_{n=2}^{\infty}(|a_n|+|b_n|)r^n(\alpha)+\frac{1}{1-r^N(\alpha)}\sum_{n=2}^{\infty}n(|a_n|+|b_n|)^2r^{2n}(\alpha)\\&= r(\alpha)-2(1-\alpha)(r(\alpha)+\ln(1-r(\alpha)))-\frac{4(1-\alpha)^2}{1-r^N(\alpha)}(r^2(\alpha)+\ln(1-r^2(\alpha)))\\&\nonumber= 1+2(1-\alpha)(\ln 2-1)\\&=d(f(0),\partial f(\mathbb{D})).
	\end{align*}
	Hence the radius $ r(\alpha) $ is the best possible. This completes the proof.	
\end{pf} 
\begin{pf}[\bf Proof of Theorem \ref{th-2.22}]
	Let $f \in \mathcal{P}^{0}_{\mathcal{H}}(\alpha)$ be given by \eqref{e-1.11a}. Using  Lemmas \ref{lem-1.8} and \ref{lem-1.9} for $ |z|=r $, we obtain 
\begin{align}\label{e-6.26} &
	r+\sum_{n=2}^{\infty}(|a_n|+|b_n|)r^n+\left(\frac{1}{1+|a_2|+|b_2|}+\frac{r^m}{1-r^m}\right)\sum_{n=3}^{\infty}n^{m-1}(|a_n|+|b_n|)^mr^{mn}\\ &\nonumber\leq r+\sum_{n=2}^{\infty}\frac{2(1-\alpha)r^{n}}{n}+\left(\frac{1}{1+|a_2|+|b_2|}+\frac{r^m}{1-r^m}\right)\sum_{n=3}^{\infty}\frac{2^m(1-\alpha)^m(r^m)^n}{n}.
\end{align}
Let
\begin{align*}
	Q_m(r):=\frac{1}{1+|a_2|+|b_2|}+\frac{r^m}{1-r^m}.
\end{align*}
A computation using $ Q_m(r) $ shows that
\begin{align}
\label{e-6.27} &
	r+\sum_{n=2}^{\infty}\frac{2(1-\alpha)r^{n}}{n}+Q_m(r)\sum_{n=3}^{\infty}\frac{2^m(1-\alpha)^m(r^m)^n}{n}\\&=\nonumber r-2(1-\alpha)(r+\ln(1-r))-2^m(1-\alpha)^mQ_m(r)\nonumber\left(r^m+\frac{r^{2m}}{2}+\ln(1-r^m)\right)\\& \nonumber\leq 1+2(1-\alpha)(\ln 2-1)
\end{align}
for $ r\leq r_m(\alpha) $, where $ r_m(\alpha) $ is the smallest root of of $ H_8(r) $ in $ (0,1) $ and 
\begin{align*}
	& H_8(r):=\\&r-1-2(1-\alpha)\left(r-1+\ln (2-2r)\right)-2^m(1-\alpha)^mQ_m(r)\left(r^m+\frac{r^{2m}}{2}+\ln(1-r^m)\right).
\end{align*}
Then we have $ H_8(r_m(\alpha))=0 $. That is
\begin{align}\label{e-6.28} &
r_m(\alpha)-1-2(1-\alpha)\left(r_m(\alpha)-1+\ln (2-2r_m(\alpha))\right)\\&\nonumber-2^m(1-\alpha)^mQ_m(r)\left(r^m_m(\alpha)+\frac{r^{2m}_m(\alpha)}{2}+\ln(1-r^m_m(\alpha))\right)=0.
\end{align}
With the help of \eqref{e-3.3}, \eqref{e-6.26} and \eqref{e-6.27} for $ |z|=r\leq r_m(\alpha) $, we obtain
\begin{align*} &
	r+\sum_{n=2}^{\infty}(|a_n|+|b_n|)r^n+\left(\frac{1}{1+|a_2|+|b_2|}+\frac{r^m}{1-r^m}\right)\sum_{n=3}^{\infty}n^{m-1}(|a_n|+|b_n|)^mr^{mn}\\&\leq d(f(0),\partial f(\mathbb{D})).
\end{align*}
To show that $ r_m(\alpha) $ is the best possible, we consider the function $ f=f_{\alpha} $ defined by \eqref{e-6.7}. For the function $ f=f_{\alpha} $ and $ |z|=r_m(\alpha) $, an simple computation using \eqref{e-3.7} and \eqref{e-6.28} shows that 
\begin{align*} &
	r_m(\alpha)+\sum_{n=2}^{\infty}(|a_n|+|b_n|)r^n_m(\alpha)+\left(\frac{1}{1+|a_2|+|b_2|}+\frac{r^m}{1-r^m_m(\alpha)}\right)\\&\quad\quad\times \sum_{n=3}^{\infty}n^{m-1}(|a_n|+|b_n|)^mr^{mn}_m(\alpha)\\&= r_m(\alpha)-2(1-\alpha)(r_m(\alpha)+\ln(1-r_m(\alpha)))\\&\quad\quad-2^m(1-\alpha)^mQ_m(r_m(\alpha))\left(r^m_m(\alpha)+\frac{r^{2m}_m(\alpha)}{2}+\ln(1-r^{m}_m(\alpha))\right)\\&\nonumber= 1+2(1-\alpha)(\ln 2-1)\\&=d(f(0),\partial f(\mathbb{D})).
\end{align*}
Therefore, the radius $ r_m(\alpha) $ is the best possible.	This completes the proof.
\end{pf} 
\begin{pf}[\bf Proof of Theorem \ref{th-2.25}]
		Let $f \in \mathcal{P}^{0}_{\mathcal{H}}(\alpha)$ be given by \eqref{e-1.11a}. By using Lemmas \ref{lem-1.8} and \ref{lem-1.9}, for $ |z|=r $, we obtain 
	\begin{align}\label{e-6.29} &
		r+\frac{[1-(1+|a_2|+|b_2|-(|a_2|+|b_2|)^2)]r}{1-(|a_2|+|b_2|)r}+\sum_{n=3}^{\infty}(|a_n|+|b_n|)r^{n}\\[2mm]&\nonumber\leq r+\frac{[1-(1+|a_2|+|b_2|-(|a_2|+|b_2|)^2)]r}{1-(|a_2|+|b_2|)r}+\sum_{n=3}^{\infty}\frac{2(1-\alpha)r^{n}}{n}.
	\end{align}
	Let $ {R}(r) $ be defined
	\begin{equation*}
		{R}(r)=\frac{[1-(1+|a_2|+|b_2|-(|a_2|+|b_2|)^2)]r}{1-(|a_2|+|b_2|)r}.
	\end{equation*}
	 A simple computation shows that
	 \begin{align}	\label{e-6.30}\;\;\;
	 	r+{R}(r)+\sum_{n=3}^{\infty}\frac{2(1-\alpha)r^{n}}{n}&= r+{R}(r)-2(1-\alpha)\left(r+\frac{r^2}{2}+\ln(1-r)\right)\\&\leq 1+2(1-\alpha)(\ln 2-1)\nonumber
	 \end{align}
	for $ r\leq r(\alpha) $, where $ r(\alpha) $ is the smallest root of $ H_9(r) $	in $ (0,1) $ and $ H_9 : [0,1)\rightarrow \mathbb{R} $ be defined by 
	$$ H_9(r):=	r+{R}(r)-2(1-\alpha)\left(r+\frac{r^2}{2}-1+\ln(2-2r)\right)-1.
	$$
 Then we have $ H_9(r(\alpha))=0 $. That is 
	\begin{align}\label{e-6.31} 
		r(\alpha)+{R}(r(\alpha))-2(1-\alpha)\left(r(\alpha)+\frac{r^2(\alpha)}{2}-1+\ln(2-2r(\alpha))\right)-1=0.
	\end{align}
	From \eqref{e-3.3}, \eqref{e-6.29} and \eqref{e-6.30} for $ |z|=r\leq r(\alpha) $, we obtain
	$$ 
	r+\frac{[1-(1+|a_2|+|b_2|-(|a_2|+|b_2|)^2)]r}{1-(|a_2|+|b_2|)r}+\sum_{n=3}^{\infty}(|a_n|+|b_n|)r^{n}\leq d(f(0),\partial f(\mathbb{D})).
	$$ 
	To show that the radius $ r(\alpha) $ is the best possible, we consider the function $ f=f_{\alpha} $ defined in \eqref{e-6.7}. For $ f=f_{\alpha} $ and $ |z|=r(\alpha) $, a simple calculation using \eqref{e-3.7} and \eqref{e-6.31} shows that 
	\begin{align*} &
		r(\alpha)+\frac{\left(1-(1+|a_2|+|b_2|-(|a_2|+|b_2|)^2)\right)r(\alpha)}{1-(|a_2|+|b_2|)r(\alpha)}+\sum_{n=3}^{\infty}(|a_n|+|b_n|)r^n(\alpha)\\&= r(\alpha)+{R}(r(\alpha))-2(1-\alpha)\left(r(\alpha)+\frac{r^2(\alpha)}{2}+\ln(1-r(\alpha))\right)\\&\nonumber= 1+2(1-\alpha)(\ln 2-1)\\&=d(f(0),\partial f(\mathbb{D})).
	\end{align*}
	Hence the radius $ r(\alpha) $ is the best possible. This comples the proof.	
\end{pf} 
 \begin{pf}[\bf Proof of Theorem \ref{th-2.31}]
 	Let $f \in \mathcal{P}^{0}_{\mathcal{H}}(\alpha)$ be given by \eqref{e-1.11a}.
 	It is not difficult to show that
 	\begin{align*}
 		|h^{\prime}(z)|\leq \alpha+(1-\alpha)\left(\frac{1+r}{1-r}\right).
 	\end{align*} 
 Therefore, we have
 \begin{align*}
 	|{J}_f(z)|=|h^{\prime}(z)|^2-|g^{\prime}(z)|^2\leq |h^{\prime}(z)|^2.
 \end{align*}
 	Using Lemmas \ref{lem-1.8} and \ref{lem-1.9}, for $ |z|=r_N(\alpha) $, we obtain 
 	\begin{align}\label{e-6.32} &
 		|f(z)|+\sqrt{|{J}_f(z)|}r+\sum_{n=N}^{\infty}(|a_n|+|b_n|)r^{n} \\&\nonumber\leq r+2(1-\alpha)(-r-\ln (1-r))+\left(\alpha+(1-\alpha)\left(\frac{1+r}{1-r}\right)\right)r\\&\nonumber\quad\quad+ 2(1-\alpha)\left(-r-\frac{r^2}{2}-\cdots-\frac{r^{N-1}}{N-1}-\ln(1-r)\right)\\&\nonumber=r-2(1-\alpha)\left(2r+\frac{r^2}{2}+\cdots+\frac{r^{N-1}}{N-1}+2\ln(1-r)\right)+\left(\alpha+(1-\alpha)\left(\frac{1+r}{1-r}\right)\right)r.
 	\end{align}
 	A simple calculations shows that
 	\begin{align}
 		\label{e-6.33} &
 		r-2(1-\alpha)\left(2r+\frac{r^2}{2}+\cdots+\frac{r^{N-1}}{N-1}+2\ln(1-r)\right)+\left(\alpha+(1-\alpha)\left(\frac{1+r}{1-r}\right)\right)r\\& \nonumber\leq 1+2(1-\alpha)(\ln 2-1)
 	\end{align}
 	for $ r\leq r(\alpha) $, where $ r(\alpha) $ is the smallest root of $ H_{10}(r) $ in $ (0,1) $. Here $ H_{10} : [0,1)\rightarrow \mathbb{R} $ defined by 
 	\begin{align*}
 		 H_{10}(r):&=r-1-2(1-\alpha)\left(2r-1+\frac{r^2}{2}+\cdots+\frac{r^{N-1}}{N-1}+\ln 2+2\ln(1-r)\right)\\&\quad\quad+\left(\alpha+(1-\alpha)\left(\frac{1+r}{1-r}\right)\right)r.
 	\end{align*}
 	Therefore, we have $ H_{10}(r_{_N}(\alpha))=0 $, which shows that 
 	\begin{align}\label{e-6.34} &
 		r_N(\alpha)-1-2(1-\alpha)\left(2r_N(\alpha)-1+\frac{r^2_N(\alpha)}{2}+\cdots+\frac{r^{N-1}_N(\alpha)}{N-1}+\ln 2+2\ln(1-r_N(\alpha))\right)\\&\nonumber\quad\quad+\left(\alpha+(1-\alpha)\left(\frac{1+r_N(\alpha)}{1-r_N(\alpha)}\right)\right)^2=0,
 	\end{align}
 	Using \eqref{e-3.3}, \eqref{e-6.32} and \eqref{e-6.33} for $ |z|=r\leq r_{_N}(\alpha) $, we obtain
 	$$ 
 	|f(z)|+|{J}_f(z)|+\sum_{n=N}^{\infty}(|a_n|+|b_n|)r^{n}\leq d(f(0),\partial f(\mathbb{D})).
 	$$ 
 	In order to show that $ r_{_N}(\alpha) $ is sharp, we consider the function $ f=f_{\alpha} $ defined by \eqref{e-6.7}. For $ f=f_{\alpha} $ and $ |z|=r_{_N}(\alpha) $, a simple calculation using \eqref{e-3.7} and \eqref{e-6.34} shows that 
 	\begin{align*} 
 		& |f(z)|+\sqrt{|{J}_f(z)|}r+\sum_{n=N}^{\infty}(|a_n|+|b_n|)r^{n}_{_N}(\alpha)\\&= r_{_N}(\alpha)-2(1-\alpha)\left(2r_{_N}(\alpha)+\frac{r^2_{_N}(\alpha)}{2}+\cdots+\frac{r^{N-1}_{_N}(\alpha)}{N-1}+2\ln(1-r_{_N}(\alpha))\right)\\&\quad\quad+\left(\alpha+(1-\alpha)\left(\frac{1+r_{_N}(\alpha)}{1-r_{_N}(\alpha)}\right)\right)r\\&\nonumber= 1+2(1-\alpha)(\ln 2-1)\\&=d(f(0),\partial f(\mathbb{D})).
 	\end{align*}
 	Therefore, the radius $ r_{_N}(\alpha) $ is the best possible.	This completes the proof.
 \end{pf}

\noindent\textbf{Acknowledgment:}  The first author is supported by the Institute Post Doctoral Fellowship of IIT Bhubaneswar, India, the second author is supported by SERB-MATRICS, and third author is supported by CSIR, India.

\end{document}